\documentclass{article}

\usepackage{amsthm}
\newtheorem{theorem}{Theorem}[section] 

\newtheorem{proposition}[theorem]{Proposition}
\newtheorem{remark}[theorem]{Remark}

\usepackage[algo2e, ruled]{algorithm2e}
\usepackage{latexsym}
\usepackage{amsmath}
\usepackage{amssymb,amsfonts,amscd}
\usepackage{graphicx,float}
\usepackage{bbding}
\usepackage{color}

\newcommand{\bC}{{\mathbb C}}

\usepackage[]{subfig}
\usepackage{stfloats}

\usepackage{epstopdf}

\usepackage{listings}
\newcommand{\DFile}[3]{{
\begin{center}
\begin{minipage}{0.9\linewidth}
\lstinputlisting[frame=single,
	breaklines=true,
	caption={#1},
	label=#2]{#3}
\end{minipage}
\end{center}
}}

\title{{\tt {P}aramotopy}: Parameter homotopies in parallel}

\usepackage{authblk}
\usepackage{hyperref}

\author[1]{
Daniel J. Bates
	\thanks{This material is based upon work supported by the National Science Foundation under Grants No. DMS-1025564 and DMS-1115668.}
}

\author[2]{
Danielle A. Brake
	\thanks{This research utilized the CSU ISTeC Cray HPC System supported by NSF Grant CNS-0923386.}
}

\author[3]{
Matthew E. Niemerg
}

\affil[1]{Colorado State University~\\~\url{bates@math.colostate.edu}}
\affil[2]{University of Wisconsin - Eau Claire \\ \url{brakeda@uwec.edu}}
\affil[3]{\url{research@matthewniemerg.com}}

\begin{document}
\maketitle

\thispagestyle{empty}
\pagestyle{empty}

\begin{abstract}
Numerical algebraic geometry provides a number of efficient tools for approximating the solutions of 
polynomial systems.  One such tool is the parameter homotopy, which can be an extremely efficient method 
to solve numerous polynomial systems that differ only in coefficients, not monomials.  This technique is frequently  
used for solving a parameterized family of polynomial systems at multiple parameter values.
Parameter homotopies have recently been useful in several areas of application 
and have been implemented in at least two software packages.  This article describes Paramotopy, a new, 
parallel, optimized implementation of this technique, making use of the 
Bertini software package.  The novel features of this implementation, not available elsewhere, 
include allowing for the simultaneous solutions of arbitrary polynomial systems in a parameterized family on 
an automatically generated (or manually provided) mesh in the parameter space of coefficients, 
front ends and back ends that are easily specialized to particular classes of problems, and 
adaptive techniques for solving polynomial systems near singular points in the parameter space.  
This last feature automates and simplifies a task that is important but often misunderstood by 
non-experts.
\end{abstract}

\section{Introduction}\label{sec:intro}

The methods of numerical algebraic geometry provide a means for approximating 
the solutions of a system of polynomials $F:\mathbb C^N\to\mathbb C^n$, i.e., those points (perhaps forming positive-dimensional components--curves, surfaces, etc.) 
$z\in\mathbb C^N$ such that $F(z)=0$.  There are many variations on these methods, but the key point is that 
polynomial systems of moderate size can be solved efficiently via homotopy continuation-based methods.  
In the case of a parameterized family of 
polynomial systems $F:\mathbb C^N\times \mathcal{P} \to \mathbb C^N$, where  the 
coefficients are polynomial in the parameters $p\in \mathcal{P}\subset \mathbb{C}^M$, a particularly efficient technique comes 
into play:  the parameter homotopy~\cite{MS89}\footnote{In fact, this technique applies when the coefficients are {\em holomorphic functions} of the 
parameters~\cite{MS89}, but we restrict to the case of polynomials as Bertini is restricted to polynomials.}.  

Parameter homotopies are quite powerful for many classes of problems that arise in practice; before detailing the theory of parameter homotopies, we first describe the basics of using standard, non-parameter homotopies.
The process of using a standard homotopy to solve a system $F$ begins with the construction of a polynomial system $G$ that is easily 
solved.  Once the system $G$ is solved, the solutions of $G$ are tracked numerically by predictor-corrector methods as the polynomials of $G$ are transformed  
into those of $F$.  Thanks to the underlying geometry, discussed for example in~\cite{SW05} or~\cite{BertiniBook}, we are guaranteed to find a superset $\widehat V$ of the 
set $V$ of isolated solutions of $F$.  The set $\widehat V$ is easily trimmed down to $V$ in a post-processing step~\cite{LDT}.

One of the two extremes among the many choices of homotopy constructions is the {\em total degree homotopy}, which typically requires
the tracking of many more paths than the number of solutions of $F$.  At the other extreme, {\em polyhedral homotopies} require the tracking of the 
exact number of paths as the number of solutions of $F$, under the assumption that the coefficients of $F$ are generic, but even for systems of moderate size, this reduction in the number of paths comes at the cost of significant computational 
effort to solve $G$.  Many possible homotopies exist between these two, each with varying levels of complexity for solving $G$ and varying numbers of paths 
to be tracked.
 
{\em Parameter homotopies} behave nicely, as the number of paths to be followed is exactly equal to the number of isolated solutions of $F(z,p)$ for almost all values 
of $p\in \mathcal P$ (under the common assumption that $\mathcal P$ has positive volume in its ambient Euclidean space) {\em and} the solution of a single $G$ will work for almost all values of $p\in\mathcal P$, so only one round of precomputation is needed regardless 
of the number of polynomial systems to be solved.  This is described in more detail in~\S\ref{s:hom}.

Parameter homotopies are not new and have been used in several areas of application~\cite{kinematics,he2013exploring,Newell,KKT} 
and implemented in at least two software packages for solving polynomial systems:  Bertini~\cite{Bertini} and 
PHCpack~\cite{PHC}.  These implementations allow the user to run a single parameter homotopy from one parameter value $p_0$ 
with known solutions to the desired parameter value, $p_1$, with the solutions at $p_0$ provided by the user.

The new software package that is the focus of this article differs from these other two implementations in the following ways:
\begin{enumerate}
\item Paramotopy accepts as input the general form of the parameterized family $F(z,p)$ ($p$ given as indeterminates), chooses 
a random $p_0\in \mathcal P$, and solves $F(z,p_0)$ via a Bertini run\footnote{Bertini provides this functionality as well.}; 
\item Paramotopy builds a mesh in the parameter space given simple instructions from the user (or uses a user-provided set 
of parameter values) and performs parameter homotopy runs from $p_0$ to each other $p$ in the mesh; 
\item Paramotopy carries out all of these runs in parallel, as available\footnote{Bertini and PHCpack both have parallel versions, 
but not for multiple parameter homotopy runs.}; 
\item Paramotopy includes adaptive schemes to automatically attempt to find the solutions of $F(z,p)$ from starting points other than 
$p_0$, if ill-conditioning causes path failure in the initial attempt; and 
\item Paramotopy is designed to simplify the creation of front ends and back ends specialized for particular applications.
\end{enumerate}

The purpose of this article is two-fold:  to provide a refresher on parameter homotopies and to describe the software package Paramotopy.
Parameter homotopies are described in more technical detail in the next section, followed by implementation details of 
Paramotopy in \S\ref{s:poly}.  Finally, a few examples and timings are provided in \S\ref{s:ex}.

\section{Homotopies}\label{s:hom}

In this section, we introduce  homotopy continuation (\S\ref{ss:homcont}), then the special setting of parameter homotopies (\S\ref{ss:paramhom}).

\subsection{Homotopy continuation}\label{ss:homcont}

Given a polynomial system $F:\mathbb C^N \to \mathbb C^N$ to be solved, standard homotopy continuation consists of three basic steps:
\begin{enumerate}
\item Choose a {\em start system} $G:\mathbb C^N \to\mathbb C^N$ similar in some way to $F(z)$ that is ``easy'' to solve;
\item Find the solutions of $G(z)$ and form the new homotopy function $H:\mathbb C^N\times \mathbb C\to \mathbb C^N$ 
given by $H(z,t)=F(z)\cdot (1-t)+G(z)\cdot t\cdot \gamma$, where $\gamma\in\mathbb C$ is randomly chosen; and
\item Using predictor-corrector methods (and various other numerical routines~\cite{AMP,AMP2,ode,SW05}), track the solutions of $G$ at 
$t=1$ to those of $F$ at $t=0$.
\end{enumerate}

There are many variations on this general theme, but we focus here on the basic ideas, leaving details and alternatives to the references.

\begin{remark}
\begin{enumerate}
\item There are simple methods for making a nonsquare polynomial system ($N\neq n$) square ($N=n$).  For simplicity in this article, we assume that 
the polynomial system is square.  See~\cite{SW05,BertiniBook} for details.
\item A discussion of the choice of an adequate start system $G$ goes beyond the scope of this paper.  It is enough to know that there are 
several such options~\cite{mhom,Li03,SW05,BertiniBook}.  We present the simplest homotopy that can be constructed, the {\em total degree homotopy}, after this remark.
\item Notice that $H(z,t)$ has the property that $H(z,1)=G(z)$ and $H(z,0)=F(z)$.
\item The extra $\gamma$ included in homotopy function $H(z,t)$ introduces randomness to the paths to be followed.  This ``gamma trick'' is discussed later in this section 
and is central to the {\em probability one} nature of homotopy continuation methods. 
\end{enumerate}
\end{remark}

As a simple example of the choice of $G$, let's consider the total degree homotopy.
Let $d_1,\ldots,d_N$ denote the degrees of the polynomials of $F(z)$.  One instance of a {\em total degree} or {\em B\'ezout} homotopy is 
\begin{equation*}
\begin{aligned}
g_1(z) &= z_1^{d_1}-1\\
g_2(z) &= z_2^{d_2}-1\\
&\ldots\\
g_N(z) &= z_N^{d_N}-1,\\
\end{aligned}
\end{equation*}
with $z=\left(z_1,z_2,\ldots,z_N\right)\in\mathbb C^N$.  This system has $d_1\cdot d_2\cdot\ldots\cdot d_N$ trivially-computed solutions, so 
a homotopy using this as a start system would have that number of paths, regardless of the number of solutions of $F(z)$.  A generalization of this sort of start system, the 
{\em multihomogeneous} or {\em $m$-hom} start system, is the standard in Bertini~\cite{Bertini} and is thus the main type of start system used 
in our implementation; see~\cite{SW05,BertiniBook} for more on that particular choice.  Regeneration~\cite{regen} is a recent advance that 
will also greatly increase the efficiency of non-parameter homotopy runs. 

\begin{figure}[H]
\begin{center}
\includegraphics[width = 3.5in]{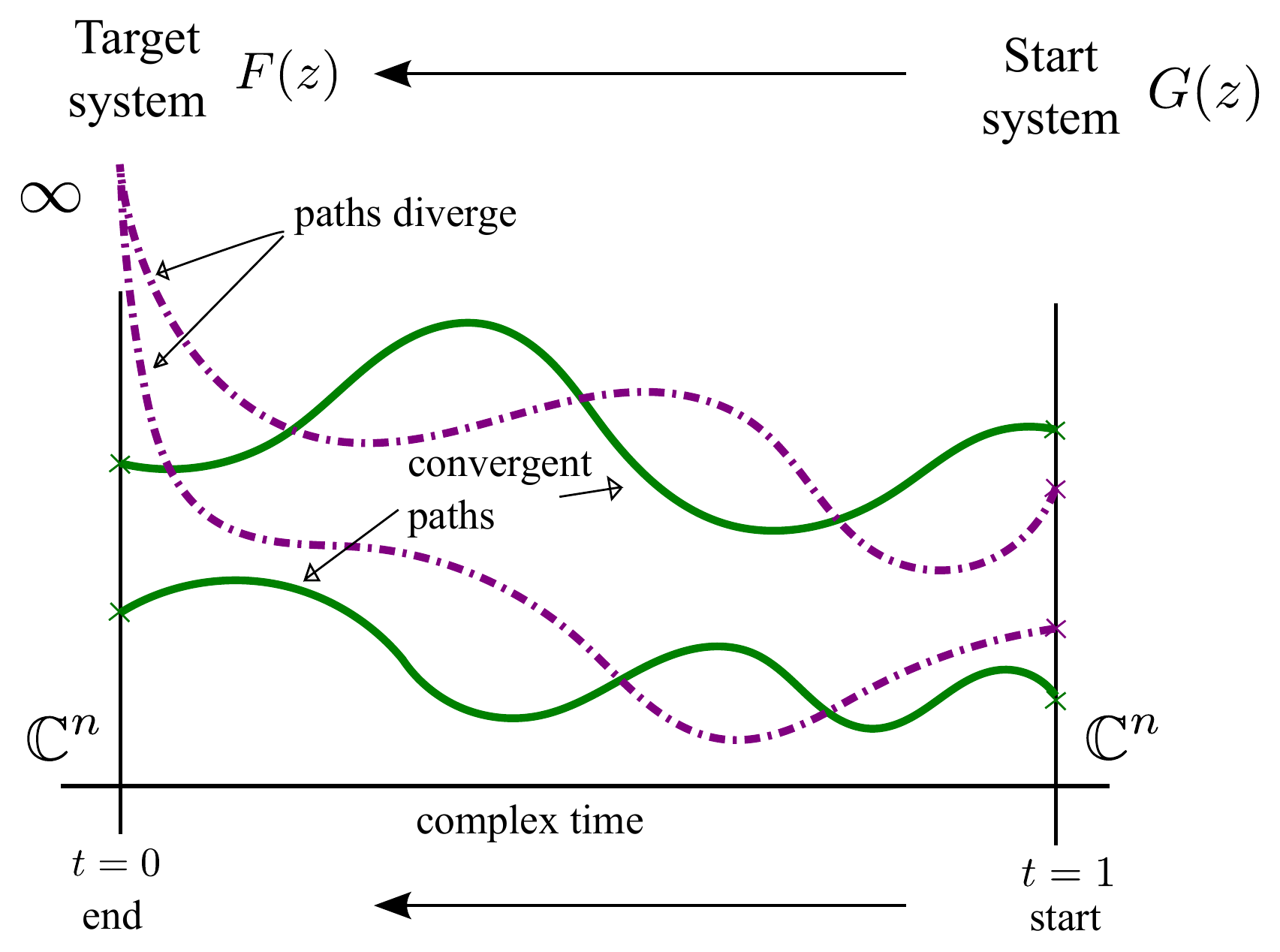}

\caption{A schematic depiction of a homotopy from system $F$ to system $G$.  There are four solutions of $G(z)$ at $t=1$.  Two solution paths 
diverge as $t\to 0$, while the other two lead to solutions of $F$ at $t=0$.}
\label{fig:homotopy}

\end{center}
\end{figure}

Once $G(z)$ is solved and $H(z,t)$ is formed, the solutions of $H(z,t)$ for varying values of $t$ may be visualized as curves.  Indeed, as $t$ 
varies continuously, the solutions of $H(z,t)$ will vary continuously, so each solution sweeps out a curve or path (also sometimes called a 
{\em solution curve} or {\em solution path}) as $t$ moves from 1 to 0.  A schematic of four such paths is given in Figure~\ref{fig:homotopy}.

Predictor-corrector methods are used to follow the solutions of $G$ to those of $F$ along these paths.  For example, 
an Euler (tangent) predictor will find a point $z^*$ near a solution of $H(z,t^*)$ for some given $t^*<1$, after which Newton's method may 
be used to correct $z^*$ back to the solution at $t^*$.  This process is repeated to move forward along the path, towards $t=0$.  
This is a vast oversimplification of what we refer to as {\em standard homotopy methods} (as opposed to parameter homotopies); indeed, in Bertini, 
Euler's method has been replaced with more accurate predictors~\cite{ode}, potential path-crossings are handled with adaptive multiprecision 
path tracking techniques~\cite{AMP,AMP2}, powerful numerical techniques called endgames are employed near the end of the 
path ($t\approx 0$)~\cite{cauchy,pseg,cauchy2}, and, for large problems, all of this is done in parallel~\cite{BertiniBook}.  

Despite all the modern safeguards against numerical difficulties, problems can still arise.  If a path is not successfully tracked 
from $t=1$ all the way to a solution at $t=0$, we refer to this as a {\em path failure}.  Causes for path failures in Paramotopy are inherited directly 
from those of Bertini~\cite{BertiniBook}, as Paramotopy relies on Bertini for all path tracking.  There are numerous such causes, though one of the most 
common comes from the situation of having two paths come near one another.  In that case, the Jacobian matrix of the homotopy function becomes 
ill-conditioned, causing an increase in precision and/or reductions to the step size.  Precision has a maximum allowed value within Paramotopy 
and Bertini, while the step size has a minimum allowed value.  If either of these thresholds is broken, the path is declared a failure.   It is important 
to note that path failures are not failures of the methods or the software.  On the contrary, a path failure is a signal to the user that the geometry 
of the path is somehow particularly tricky and that more care needs to be taken.   The mitigation of common path failures for parameter 
homotopies is described in~\S\ref{s:step2fails}.

\subsection{Parameter homotopies}\label{ss:paramhom}

Suppose we wish to solve a parameterized polynomial system $F(z,p)$ in variables $z$ and parameters $p$ at a (possibly very large) number of 
points in parameter space, i.e., we want to find $z$ such that $F(z,p')=0$ for varying values $p=p'$.  If we 
know all isolated, finite, 
complex solutions at some generic point $p=p_0$ in a {\em convex}\footnote{Handling non-convex parameter spaces is significantly more difficult and is described later.} 
parameter space $\mathcal P$, the underlying theory allows us to make use of a {\em parameter} or 
{\em coefficient-parameter homotopy}~\cite{MS89}.
The usefulness of this software becomes readily apparent from the following proposition, proved in somewhat different language in~\cite{SW05}.  The proposition guarantees that we can find 
the isolated, finite, complex solutions of $F(z,p')$ simply by following paths through the parameter space, $\mathcal P \subset \bC^M$, from the solutions of $F(z,p_0)$.  

\begin{proposition}
The number of finite, isolated solutions of $F(z,p)$ is the same for all $p\in\mathcal P$ except for a measure zero, algebraic subset $\mathcal B$ 
of $\mathcal P$.
\label{prop}
\end{proposition}

This proposition gives us a probability one guarantee that a randomly chosen path 
through parameter space will avoid $\mathcal B$.  
Furthermore, assuming $\mathcal P$ is convex, a straight line segment through parameter space from a randomly chosen 
$p_0\in\mathcal P$ to a prespecified target $p_1\in\mathcal P$ will, with probability one, not pass through the set $\mathcal B$.  
 By moving in a straight line from a random starting point in parameter space, we should not have 
any path-crossings or divergent paths, i.e., the straight line from our starting point to $t=0$ should miss $\mathcal B$.

This immediately implies a (known) technique for solving many polynomial systems from the same parameterized family with parameter 
space~$\mathcal P$.  First, find all finite, isolated, complex solutions for some randomly chosen $p_0\in\mathcal P$.  We refer 
to this as Step 1.  Second, for each parameter value of interest, $p_i\in\mathcal P$, simply follow the finite, isolated, complex solutions 
through the simple homotopy $H(z,t) = F(z,p_0)\cdot t + F(z,p_i)\cdot (1-t)$.  We refer to this as Step 2.  Notice that the randomly 
chosen $\gamma$ from standard homotopies can be neglected in this homotopy since 
$p_0$ is chosen randomly.  We describe in~\S\ref{s:step2fails} how we monitor these Step 2 runs in case paths fail and also how we 
handle such failures.

\begin{figure}[H]
\begin{center}
\includegraphics[width = 5in]{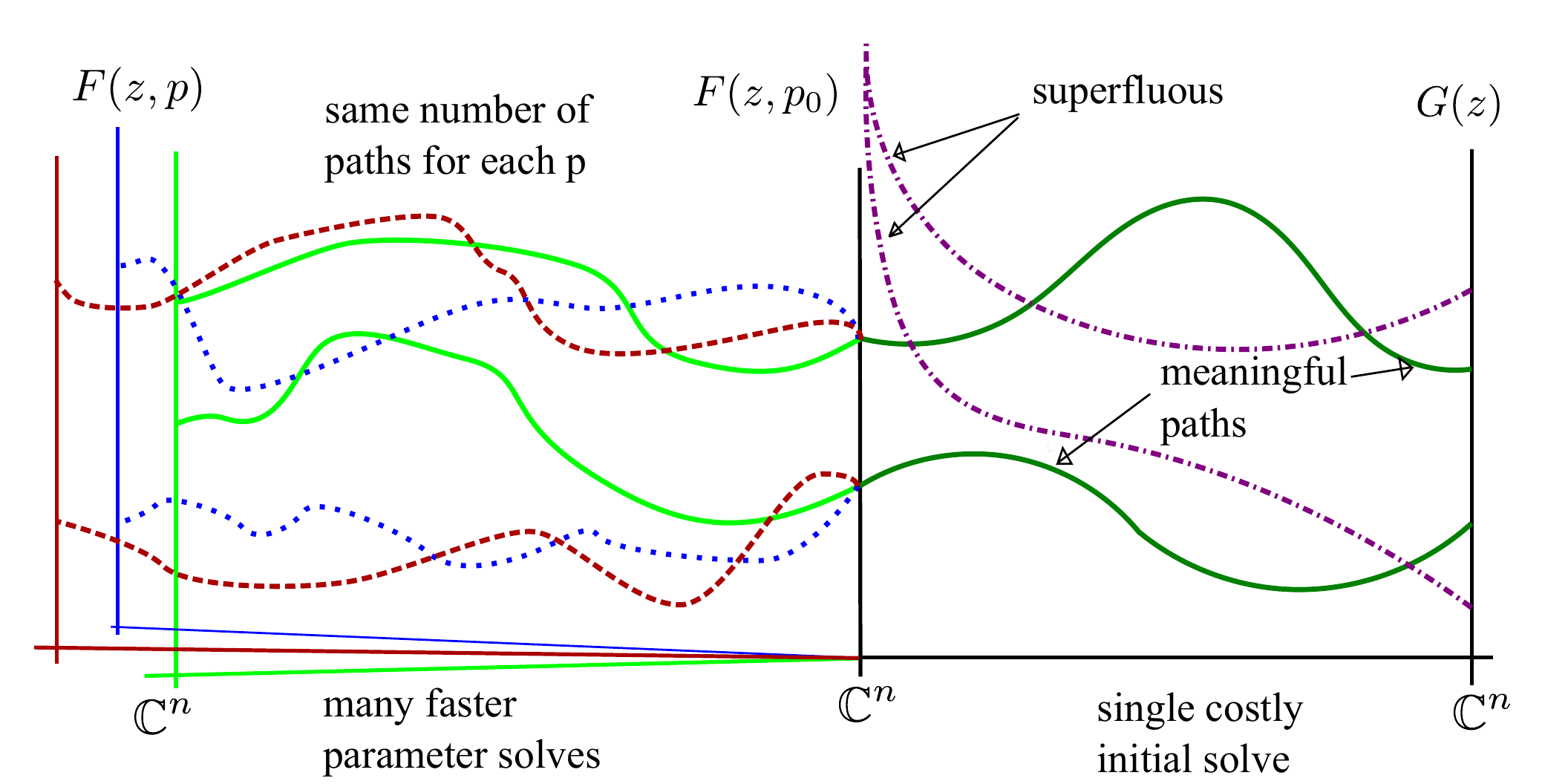}

\caption{A schematic of a parameter homotopy.  The right half of the figure corresponds to a possibly expensive Step 1 run.  The left half illustrates three parameter homotopy 
runs with two paths and two solutions each, at three different points in parameter space.}
\label{fig:parameterhomotopy}

\end{center}
\end{figure}

For the cost of a single Step 1 solve at some random point $p_0$ in the parameter space, we may rapidly solve many other polynomial 
systems in the same parameterized family.  Indeed, there are a minimal number of paths to follow in each Step 2 run and no pre-computation
cost exists beyond the initial solve at $p_0\in\mathcal P$.  

\subsubsection*{The randomness of $p_0$}  It is important to choose a random complex $p_0$ in $\mathcal P$. 
The measure zero 
set $\mathcal B$ described above contains all parameter values for 
which the number of solutions is less than the generic number of solutions.  If we were to choose $p_0$ from $\mathcal B$, we would not find 
all solutions for the parameter values of interest since we would not have as many paths as the number of solutions at each parameter value 
of interest.  By choosing $p_0$ at random, there is a zero probability that $p_0\in \mathcal B$.  If we were to use a special $p_0$, 
there is no guarantee that $p_0\notin \mathcal B$.  

\subsubsection*{The value of parameter homotopies}To get a sense of the savings from using parameter homotopies over repeated standard homotopies, suppose you can solve a parameterized 
polynomial system $F(z,p)$ for a single point $p$ in parameter space with $m$ paths.  Suppose further that you wish to solve $F(z,p)$ for 
$k$ different points in parameter space, each generically having $\ell$ solutions.  With the repeated use of a standard homotopy, you would 
need to follow a total of $km$ paths.  With a single Step 1 run at a randomly chosen complex point in the parameter space, followed by $k$ Step 
2 runs to the $k$ points of interest, you would need to follow a total of $m+k\ell$ paths.
The savings are clearly significant if $\ell << m$, especially when $k$ is large.  For example, if 
$k=1000$, $m=10000$, but $\ell=10$ (these are not exaggerated numbers), the number of paths to be tracked is reduced from 10 million to 
20 thousand by using parameter homotopies, but there is no degradation in the value of the output.  If $\ell\approx m$, it is unclear whether there 
is value in using a parameter homotopy.  Indeed, it has been noticed that paths behave differently for different types of homotopies.  This is an 
interesting open problem but is beyond the scope of this paper.

\subsubsection*{Handling non-convex parameter spaces} The above restriction that 
$\mathcal P$ be convex simplifies the discussion but is not theoretically necessary.  In our experience, it seems that {\em most} parameter spaces are 
convex, though it sometimes happens that parameter spaces may not be convex and it could very well be that, for certain application areas, parameter 
spaces are typically non-convex, perhaps not even path-connected.  

If $\mathcal P$ is not convex, it can happen that a line segment from $p_0\in\mathcal P$ to 
$p_1\in\mathcal P$ passes out of $\mathcal P$.  In that case, it could happen that the root count in the ambient Euclidean space $\mathcal X$ containing 
$\mathcal P$ is higher than that of $\mathcal P$, resulting in the failure of paths when passing from $\mathcal X\setminus \mathcal P$ back into 
$\mathcal P$. 

One potential mitigation is to replace $\mathcal P$ with $\mathcal X$, accepting that there will be more paths to follow from a generic 
$p_0\in \mathcal X  \setminus \mathcal P$.  In this case, care must be taken that paths  remain within $\mathcal P$ once they first enter $\mathcal P$, since the first entry into $\mathcal P$ could cause path failures, leading us back to the situation of the previous paragraph.  Another mitigation for 
non-convex parameter spaces is to force the path to stay within $\mathcal P$, either by parameterizing some curve from $p_0$ to $p_1$ or by choosing 
a piecewise linear path that stays within $\mathcal P$ (if one exists).  See {\em nested parameter spaces} in~\cite{BertiniBook} for more on this situation.  

In any case, it is important to note that 
\begin{enumerate}
\item Paramotopy makes the assumption that the parameter space is convex, but 
\item With care, users can handle other sorts of parameter spaces (including non-convex parameter spaces) via Paramotopy.
\end{enumerate}

In the next section, we describe the main algorithm for Paramotopy, involving one Step 1 run, followed by many Step 2 runs.  We also 
point out some of the more technical aspects of Paramotopy, such as the use of parallelization and data management, as well as the use 
of multiple parameter homotopy runs from different starting points $p_0$ to find the solutions for particular parameter values for which there were 
path failures in the initial parameter homotopy. 

Before moving on, we note that the concept of a parameter homotopy is not new.  In fact, this powerful idea dates back to the late 1980s~\cite{MS89}, with a 
somewhat more restrictive form described in~\cite{Yorke1989}.  
A thorough exposition may be found in~\cite{SW05}.

\section{Implementation details}\label{s:poly}

Paramotopy is a C++ implementation of parameter homotopies, relying heavily on Bertini~\cite{Bertini}.  In this section, 
we provide the main mathematical algorithm, Algorithm~\ref{alg:basic}, pseudocode for the fundamental algorithm, discuss how path failures 
are managed automatically, provide the 
technical details on both parallelization and data management in Paramotopy, and provide details on how to interface with 
Paramotopy from {\tt Matlab}.

\subsection{Main algorithm}\label{s:algo}

We first present the main parameter homotopy algorithm that is implemented in Paramotopy.  Note in particular the 
input value $K$ and the while loop at the end, both included to help manage path failures during the Step 2 runs.  Also, note that 
this algorithm assumes that $\mathcal P = \mathbb C^M$, for some $M$.  The use of Paramotopy for other parameter spaces is 
described in~\S\ref{s:other_param}.

\begin{algorithm2e}[t]
\SetAlgoNoLine

\KwIn{$F(z;p)$, a set of polynomial equations in $N$ variables $z\in\mathbb C^N$ and $M$ parameters $p\in L \subset \mathcal P = \mathbb C^M$; 
$\ell\ =\ \mid L\mid$ parameter values at which the solutions of $F(z;p)$ are desired; bound $K$ on the number of 
times to try to find solutions for any given $p\in L$, in the case of path failures.}

\KwOut{List of solutions of $F(z;p)=0$ for each $p\in L$.}
 Choose random $p_0\in \mathcal P $\;
 Solve $F(z;p_0)=0$ with any standard homotopy.  (Step 1)\;
 Store all nonsingular finite solutions in set $S$\;
 Set $\mathcal F:=\varnothing$.   (Beginning of Step 2.)\;  
\For{$i$=1 to $\ell$}{
  Construct parameter homotopy from $F(z;p_0)$ to $F(z;p_i)$\;
  Track all $|S|$ paths starting from points in $S$\;
  Set $\mathcal F := \mathcal F\cup \left\{ i\right\}$ if any path fails\;
} (End of Step 2.)\;
 Set $k:=0$.  (Beginning of path failure mitigation.)\;
\While{$|\mathcal F|>0$ and $k<K$} 
{
  Set $\mathcal F'=\varnothing$\;
 Choose random $p'\in \mathcal P$\;
  Solve $F(z;p')=0$ with a parameter homotopy from $p_0$\;
  \For{$m$=1 to $|\mathcal F|$}{
     Solve $F(z;p_{\mathcal F\left[m\right]})=0$ with a parameter homotopy from $p'$ to $p_{\mathcal F\left[m\right]}$\;
     Set $\mathcal F' := \mathcal F'\cup \left\{ m\right\}$ if any path fails\;
   } 
   Set $\mathcal F:=\mathcal F'$ and increment $k$\;
}  (End of path failure mitigation.)
\caption{Main algorithm.}
\label{alg:basic}
\end{algorithm2e}

\begin{remark}
To find all solutions for all $p\in L$, we must have that all solutions of $F(z,p_0)$ are nonsingular as we can only follow paths starting from nonsingular solutions during the parameter homotopies after the first run.  
Deflation~\cite{def1,def2} could be used to regularize singularities in Step 1 before beginning Step 2, but this is not currently implemented.
\end{remark}

\begin{algorithm2e}[H]
\SetAlgoNoLine

 Load input file and user preferences (otherwise use default preferences)\;
 User modifies settings for run\;
 Write \texttt{config} and \texttt{input} sections of Bertini input file  (Step 1)\;
 \texttt{system()} call to Bertini, in parallel if desired (Step 1)\;
User ensures quality of Step 1 results (see remark below)\;
Paramotopy calls \texttt{system(`mpilauncher step2')}\;
Load polynomial system information\;
Load runtime preferences\;
\If{id==0}{
	\texttt{head()}  (Algorithm~\ref{alg:head})\; }
\Else{
		\texttt{worker()}   (Algorithm~\ref{alg:worker})\;}
Write timing data to disk\;
Return to Paramotopy\;

 \texttt{FailedPathAnalysis()} \;

\caption{Paramotopy Implementation of Main Algorithm.  
\label{alg:Paramotopy}} 
\end{algorithm2e}

\begin{remark}\label{rem:param}
Paramotopy does not currently check the quality of the Step 1 results.
While theory dictates that some paths converge while others diverge as 
$t\to 0$, the reality is that paths can fail for numerical reasons.  For example, 
the path tracker can jump between paths if two paths come very near.  
This is largely mitigated in Bertini via adaptive precision and a check at $t=t_{endgame}$ 
(with $t_{endgame} = 0.1$ as the default in Bertini) 
whether all paths are still distinct.  But there is no known way to remove such 
crossings with certainty.  Thus, the user should consider reviewing the output 
of Step 1 before launching Step 2.  If nothing else, a rerun of Step 1 could 
increase confidence in the results.
\end{remark}

\begin{algorithm2e}[H]
\SetAlgoNoLine
Create input file in memory\;
Get \texttt{start} points from completed Step 1\;
MPI\_Bcast Bertini {\tt input} file and {\tt start} file to each worker\;
Get first set of points, from mesh or user-provided file\;
Distribute first round of points\;
\While{points left to run}{
	MPI\_Receive from any source\;
	 Create next set of parameters\;
	MPI\_Send new parameters to \texttt{worker()}  (Algorithm~\ref{alg:worker})\;
} 
\For {each worker}{
	 MPI\_Send kill tag\;
} 
Write timing data to disk\;
Delete temp files\;
\caption{Paramotopy head function.\index{head function}} 
\label{alg:head}
\end{algorithm2e}

\begin{algorithm2e}[H]
\SetAlgoNoLine

MPI\_Receive   {\tt start} and {\tt input} files from head\;
Change directory to working directory\;
Call initial Bertini using Start and Input, to seed memory structures\;
Collect all \texttt{.out} files from initial Bertini call\;
Receive initial round of work\;
Write \texttt{.out} files\;
\While{1}{
MPI\_Receive from head for either terminate tag, or number of points to solve\;
\If {Terminate}{done}
\Else{
	\For {Points to solve}{
		Write \texttt{num.out} file for straight line program\;
		Call Bertini library\;
		Read selected output files into buffer\;
		\If {bufferthreshold exceeded}{
			Write corresponding file buffer to disk\;
			Clear buffer\;
		} 
	} 
	MPI\_Send to head that work is completed\;
	}
} 

\caption{Paramotopy worker function.\index{worker function}} 
\label{alg:worker}
\end{algorithm2e}

\subsection{Handling path failures during Step 2}\label{s:step2fails}

If a path fails during a Step 2 run for some parameter value $p\in L$, Paramotopy will automatically 
attempt to find the solutions at $p$ by tracking from a different randomly chosen parameter value 
$p'\neq p_0\in \mathcal P$.  It will repeat this process $K$ times, with $K$ specified by the user.  
This is the content of the while loop at the end of the Main Algorithm.

The idea behind this is that paths often fail for one of two reasons: 
\begin{enumerate}
\item the path seems to be diverging, or
\item the Jacobian matrix becomes so ill-conditioned that either the steplength drops below the 
minimum allowed or the precision needed rises above the maximum allowed.
\end{enumerate}

For parameter homotopies, a path failure of the first type is possible for either of two reasons:  either 
the path really is diverging or the norm of the solution is above a particular threshold.  
In the former case, it can happen that the nature of the solution set at target value $p$ differs from that at 
a generic point in the parameter space, e.g., there could be fewer finite solutions at $p$. 
Such path failures are captured and reported by Paramotopy, but there is simply no hope for 
``fixing'' them as this result is a natural consequence of the geometry of the solution set, i.e., 
$p$ is inherently different from other points in parameter space, so Paramotopy takes the 
correct action in reporting it.  In the latter case, it can happen that  the scaling of 
the problem results in solutions that are large in some norm, e.g., $|z|_\infty > 10^5$ as is the default in the current 
version of Bertini.  If this is suspected, the user could rescale the system or adjust the threshold 
\verb1MaxNorm1 and run the problem again. 

For the second type of path failure, the ill-conditioning is caused by the presence of a singularity  
$b\in \mathcal B$ near or on the path between $p_0$ and $p$.  By choosing new starting point $p'$ ``adequately far'' from 
$p_0$, it should be feasible to avoid the ill-conditioned zone around $b$ unless $b$ is near the target value $p$.
In this last case, it is unlikely that choosing different starting points $p'$ will have any 
value, which is why we have capped the number of new starting points allowed at $K$.

For now, the new point $p'$ is chosen randomly in the unit hypercube.  Future work will detect where in parameter 
space the failures have occurred and bound $p'$ away from this region. 

Since it cannot easily be determined which paths from $p'$ to $p$ correspond to the failed 
paths from $p_0$ to $p$, there is no choice but to follow all paths from $p'$ to $p$.  To find 
all solutions at $p'$, we simply use a parameter homotopy to move the solutions at 
$p_0$ to those at $p'$.  Of course, if there are path failures, we must choose yet another $p'$ and 
try again.

\subsection{Handling parameter spaces other than $\mathbb C^M$}\label{s:other_param}

As described near the end of~\S\ref{ss:paramhom}, Paramotopy may be used to handle parameter spaces other 
than the simplest parameter space, $\mathbb C^M$ for some $M$.  However, some changes are needed in the 
algorithm.

If $\mathcal P\subset \mathbb C^M$ is a proper, convex subset of $\mathbb C^M$, Algorithm~\ref{alg:basic} needs 
only one change:  $p_0$ must be somehow chosen within $\mathcal P$.  To accommodate this, Paramotopy 
allows the user to specify $p_0$.

If $\mathcal P$ is a proper, non-convex set, more work is required.  The Step 1 run would be the responsibility 
of the user, as in the previous paragraph, and it would be up to the user to string together subsequent Paramotopy 
runs to stay within $\mathcal P$.  As this case appears to be both uncommon and highly complicated, we leave 
handling such situations the responsibility of the user.

\subsection{Parallelization and data management}

One of the features of Paramotopy that sets it apart from Bertini is the use of parallel computing for multiple 
parameter homotopies.  Bertini includes parallel capabilities for a single homotopy run, but not for a sequence of runs.
Parallelization was achieved using the head-worker paradigm, implemented with MPI.  A single process controls the distribution of parameter points to the workers, which constitute the remainder of the processes.  Workers are responsible for writing the necessary files for Bertini and for writing their own data to disk.

Bertini creates structures in memory  by parsing an {\tt input} file.  As {\tt input} is interpreted, several other files are created.  These contain the straight line program, coefficient values, variable names, etc.   Since the monomial structure of the polynomials in each Step 2 run is the same, almost all of these files are identical from one run to the next, so almost all this parsing is unnecessary.  The only file that needs to be changed between runs is the file containing parameter values.  Parsing requires a significant amount of time especially when compared to the short time needed for parameter homotopy runs, and since we call Bertini repeatedly, we eliminate as much of this parsing as possible.  We do so by calling certain Bertini functions from a compiled library, so as to prevent both the repeated parsing of an {\tt input} file and to preserve the necessary structure of the temporary files.

As Bertini runs through the paths of a homotopy, it records path-tracking data files.  To prevent proliferation in the number of files needed to contain the data from the Paramotopy run, the Bertini output data is read back into memory, and dumped into a collective data file.  The collective data files have a maximum buffer size (the default of which is 64MB), and once the buffer size is reached, the data in the buffer is written to the file, and the process repeats by storing the Bertini output data in memory until the buffer is full once more.  On modern systems, in principle all data could be collected into a single file, but transfer of data out of a computing cluster can be cumbersome with extremely large files.

Repeated writing and reading is taxing on hard drives and clogs a LAN if the workers are using network drives.  To free workers from having to physically write temporary files to electronic media storage, an option is provided to the user to exploit a shared memory location (or ramdisk), should it be available.  The default location for this is \texttt{/dev/shm}.  This is commonly available on Linux installations such as CentOS.  

As a piece of scientific software, knowledge of efficiency and performance are important.  To this end, we have developed a custom class for timing statements that records information for each parameter point by utilizing the \texttt{chrono} standard library.  This enables analysis of performance by process type, which appears below in the demonstration sections (\ref{s:ex}).  If the user is not interested in timing, the program may be compiled without the timing statements by making use of the appropriate compiler flags.

\subsection{Front ends and back ends}


Real-world problems may involve many parameters.  This could be problematic when one wants to discretize a parameter space into a uniform sample as the number of parameter points of interest can easily reach into the millions or even billions.  Hence, Paramotopy contains support for both linear uniform meshes of parameters, as well as user-defined sets of parameter values stored in a text file.  Systems with few parameters can make excellent use of computer-generated discretizations.  In contrast, systems with many parameters perhaps could use the Monte Carlo sampling method to collect useful information.  To use this functionality, the user must generate a text file containing whitespace-separated real and imaginary pairs for each parameter, with distinct parameter points being on separate lines.  An example with a user-defined parameter sampling is given in Section~\ref{ex:robo}, while computer generated regular meshes are used in Sections~\ref{ex:monks}-\ref{ex:cube}.

A generic {\tt Matlab} interface for gathering, saving, and plotting data from an arbitrary Paramotopy run is provided on the Paramotopy website.  It can handle both the mesh-style parameter discretizations, as well as user-defined sets.  The mesh-style runs generate an $n$-dimensional array, with the number of solutions, number of real solutions, and solution values stored accordingly.  Because the user-defined runs may not emit such a convenient method of storage, a 1D array is created, with line number from the parameter file corresponding to the index in the {\tt Matlab} data structure.  

Regarding plotting in this {\tt Matlab} function, two and three dimensional data sets are plotted automatically.  Higher dimension data can be plotted with the user's choice of variables appearing on the axes, with a variable used for coloring the points if desired, for display of up to 4 dimensions.  Additional basic display techniques included are movie-making and display of the number of real solutions in parameter space.  Details on how to use these features may be found in the Paramotopy user's manual, available 
from the Paramotopy website.

Other software packages, including {\tt surfex} and {\tt surfer}, can be used for viewing algebraic surfaces given by one polynomial in three 
variables~\cite{surfex,surfer}.  {\tt Bertini\_real} is another such software package that has no such restrictions on the number of polynomials and 
variables, though it will find only the real solutions within complex curves and surfaces~\cite{BertiniReal}.

\section{Examples}\label{s:ex}

In this section, we walk through a few examples of Paramotopy runs demonstrating various features.

\subsection{Kinematics}
\label{ex:robo}

Kinematics problems can often be cast in polynomial systems language; see {\em e.g.} \cite{robot_homotopy,kinematics} and a slew of others.  Robotic arms, whether prismatic or revolute in nature, can be described in terms of now-canonical parameters \cite{dh}.  Using homogeneous matrices, inverse kinematic equations can be derived for any robot, and in the case of revolute joints, are expressed in terms of sines and cosines of the joint variables.  Finally, by coupling such pairs through a Pythagorean identity, the trigonometric functions can be made polynomial.

Here, we present a very basic manipulator, purely for demonstration purposes.  Consider a three-link spatial manipulator, with each link having joint length unity.  For such a robot, Denavit Hartenberg (DH) parameters are given in Table~\ref{tab:111dh}.  To get a basic mapping of the workspace, we can sample $(x,\,y,\,z)$ coordinates randomly, and feed Paramotopy this sample.  See Input~\ref{robo}.  Running Paramotopy will then give us a map of the space, in terms of the number of real configurations possible for each point in space, as well as allowing us to estimate the volume of the workspace.

\begin{table}
\begin{center}
\caption{DH Parameters for 111 demo robot.}
\label{tab:111dh}
\begin{tabular}{ c  c  c  c }
\hline
$\theta_i$			&	$\alpha_i$			&	$a_i$ (m)		&	$d_i$ (m)		\\ \hline  \hline
$\theta_1$ 		&	$\pi/2$				& 	1			&	0		\\
$\theta_2$ 		&	0				& 	1			&	0	\\
$\theta_3$		&	$\pi/2$					&	1			&	0			\\
\hline
\end{tabular}
\end{center}
\end{table}

\begin{figure}[th]
\begin{center}
\includegraphics[width =0.56\linewidth]{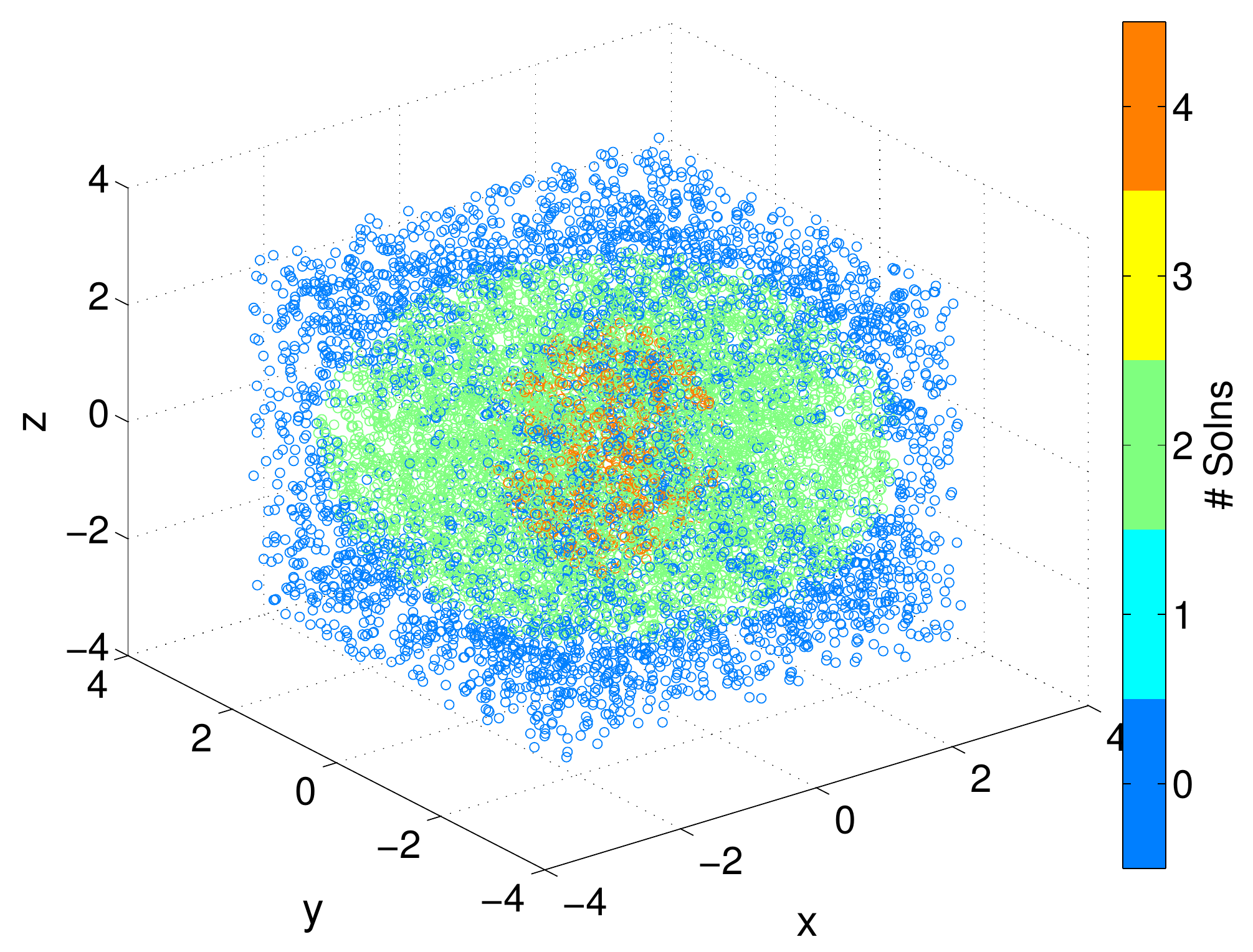}
\caption{Scatter plot of the number of real solutions versus sample point in parameter space.}
\label{fig:111nsolns_full}
\end{center}
\end{figure}

A plot of a random $10^4$ point sampling with coloring according to the number of real solutions is given in Figure~\ref{fig:111nsolns_full}.  Every point in the sample has an even number of real solutions, with a core around the origin having four, and the remainder having two.  A total of 58\% of samples had a positive number of real solutions; therefore, we can estimate that the volume of the workspace for this robot is $\approx126.6$ cubic units.  

The ability to handle arbitrary parameter samples makes Paramotopy a powerful tool.  One can perform initial sampling of a space, interpret the results using the Matlab codes provided, and make a new run using a refined sample.  Alternatively, one can perform non-linear parameter scans; {\em e.g.} one could run logarithmic samples, or sample directly on interesting sets.

We note that some level of automation is possible, as well, in that one could feed a list of numeric inputs to Paramotopy via the {\tt <} directive in the command line.  This would allow the user to pass pre-defined commands to the program, and would, for example, after the completion of a successful {\em ab ignition} run, proceed onto step 2, without the need for constant monitoring of the program.

\subsection{Dynamical Systems}
\label{ex:monks}
As derived by Ken Monks in an unpublished work, we have the `Monks Equations', which describe the amplitude of interacting waves on an annulus and are related to growth patterns in cacti and other plants.  They are interesting because the number of real solutions depends on the parameter values and because some equilibrium solutions are stable while others are unstable.  The differential equations are the following square coupled four-complex-dimensional system:
\begin{align}
\dot{z_0} &=\mu_0 z_0 + \bar{z_1} z_2 - \gamma z_0(S - |z_0|^2)\,, \nonumber \\
\dot{z_1} &=\mu_1 z_1 + \bar{z_0} z_2 + \bar{z_2} z_3 - \gamma z_1(S - |z_1|^2)\,, \label{eqn:monks}\\
\dot{z_2} &=\mu_1 z_2 + \bar{z_0} z_1 + \bar{z_1} z_3 - \gamma z_2(S - |z_2|^2)\,, \nonumber \\
\dot{z_3} &=\mu_0 z_3 + z_1 z_2 - \gamma z_3(S - |z_3|^2)\,,\nonumber 
\end{align}
where $S = 2 \sum_{j=0}^4{|z_j|^2}.$  Monks notes that the three parameters $\mu_0, \mu_1, \gamma$ must be real, and that $\gamma>0$.  Note that the complex conjugate operator is nondifferentiable and nonalgebraic.  To get around this difficulty, we seek pure real solutions ($z_i = \bar{z}_i$), and drop the bar from each equation.  See Input~\ref{monks}.  

To get a sense of the raw speedup offered by the Paramotopy method, we did the following comparison.  Using the random parameter values $p_0$ as a comparison point, and running {\tt Bertini 1.3.1} in serial mode on the same machine on which the timing runs were performed,  1000 iterations at $p_0$ took 89.4 seconds, or about 0.089 seconds per execution.  In Paramotopy, 72 workers and a single head node ran 110592 points in 221.75 seconds, including network communication and data collection, about  0.002
seconds per parameter point.  The inferred raw ratio for speedup was  44.6.  However, the parameter points in the sampling would have varying distances from the singularities, and perhaps run quicker or slower than the random point used in this timing demo, which due to genericity would run fairly quickly.  Regardless, the speedup in this example is significant, especially given the impracticality of running Bertini in parallel mode for this small problem.

\begin{figure}[th]
\begin{center}

\subfloat{\includegraphics[width =0.56\linewidth]{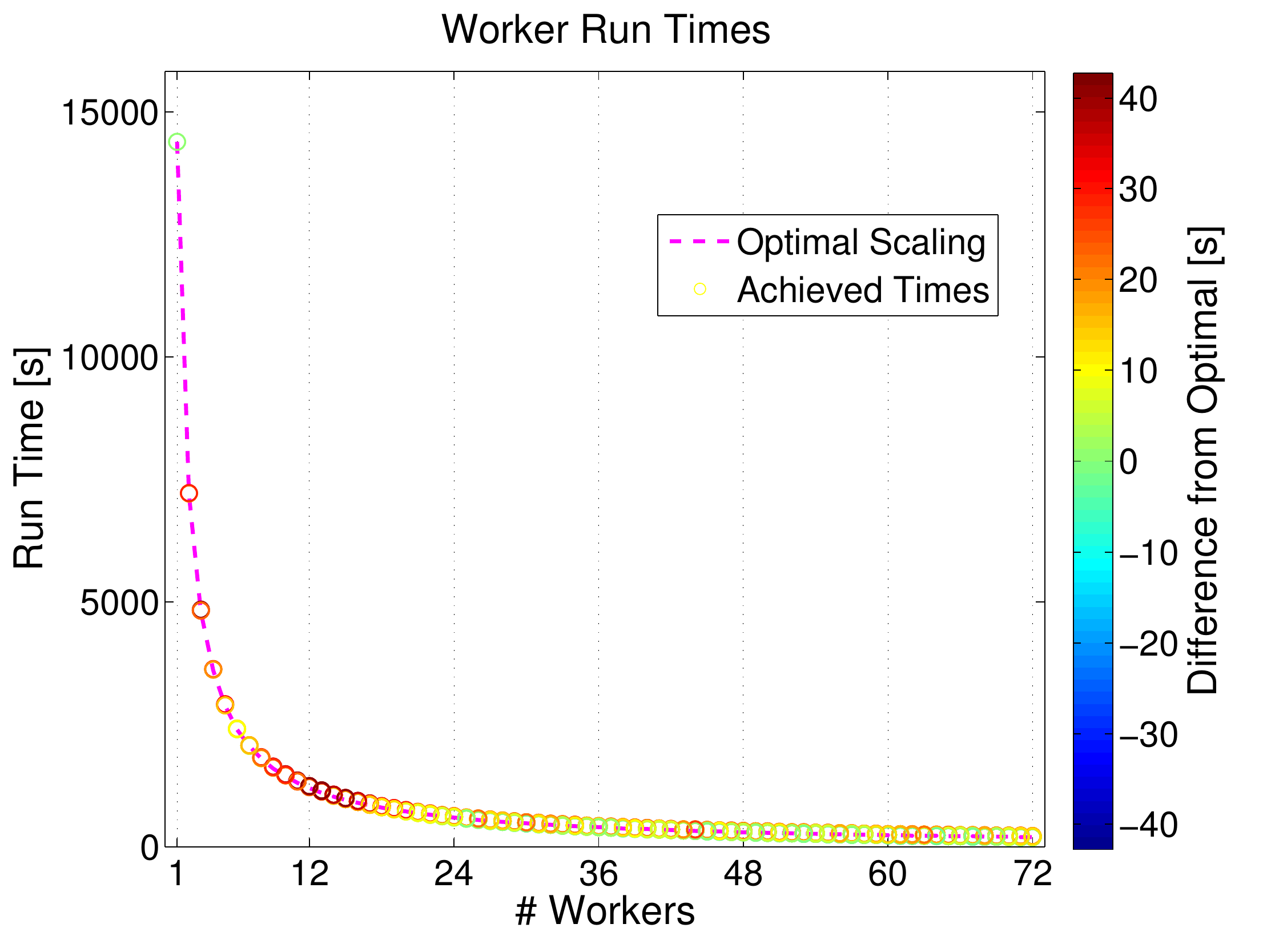}}
\subfloat{\includegraphics[width =0.46\linewidth]{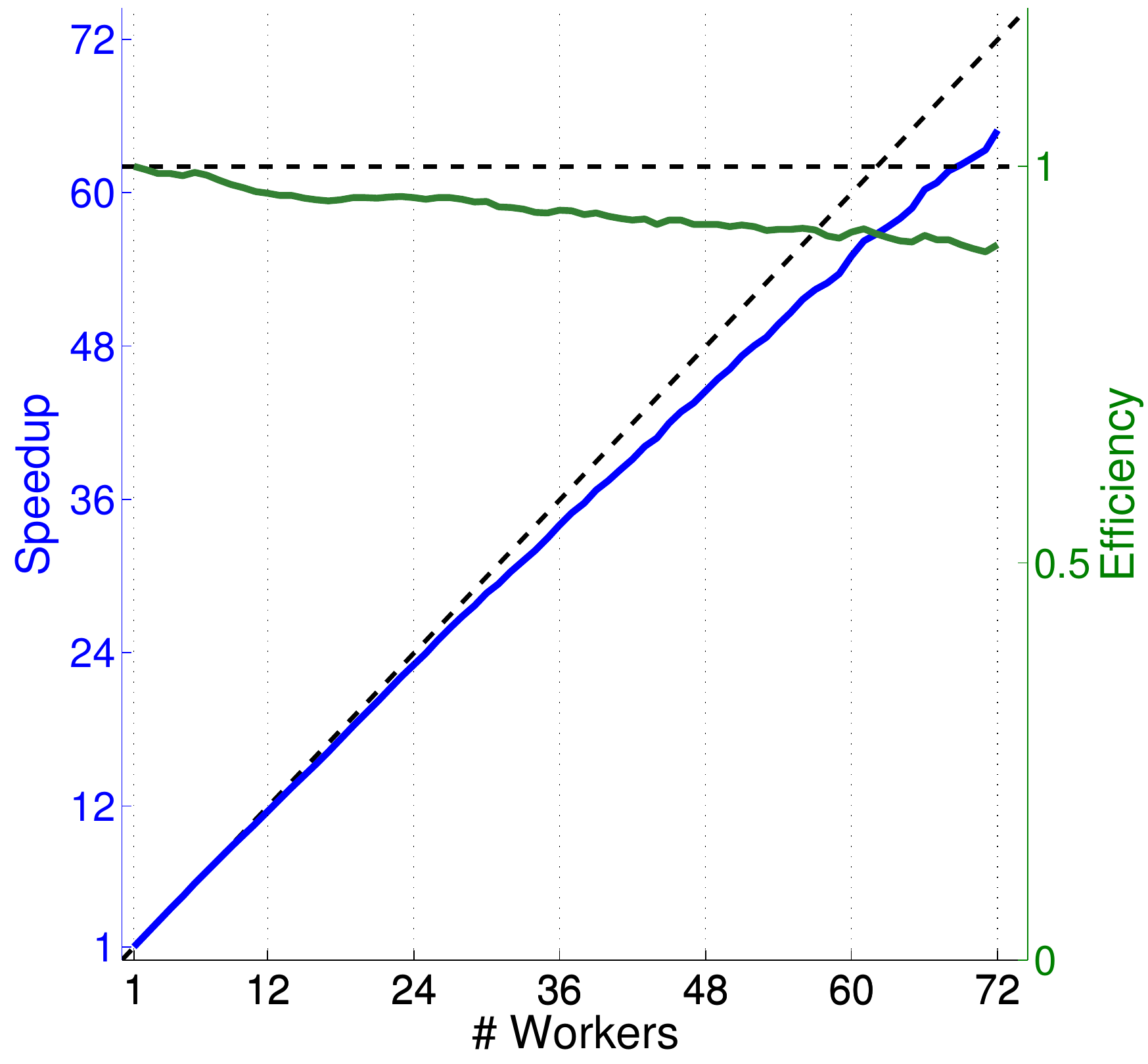}}
\caption[Paramotopy timing -- Monks Equations ]{Monks equations timing, 110592 parameter points, scaling from 1 to 72 workers. Maximum speedup of $\approx$63x achieved.}
\label{fig:monkstiming}
\end{center}
\end{figure}

Scaling tests of Paramotopy were performed using the Monks equations, on a mesh of 48 points in each of three parameters, for a total of 110592 parameter points.  Using the timing results from one worker as a base for calculation, we ran the same parameter sampling using identical $p_0$, with up to 72 workers.  Figure~\ref{fig:monkstiming} shows the results.  We saw great linear scaling throughout the tested range, with speedup and efficiency dropping off slightly as the number of workers approached 72.  The trajectory of the two curves indicates that Paramotopy would have scaled to well over 72 processors.  

We note here that scaling results depend on the problem size.  In general, the larger the parameter sampling, the more processors one can be use.  To be able to handle larger problems, Paramotopy has tunable settings for buffering memory in data collection, and parameter point distribution.

\begin{figure}[t]
\begin{center}
\includegraphics[width = 0.9\linewidth]{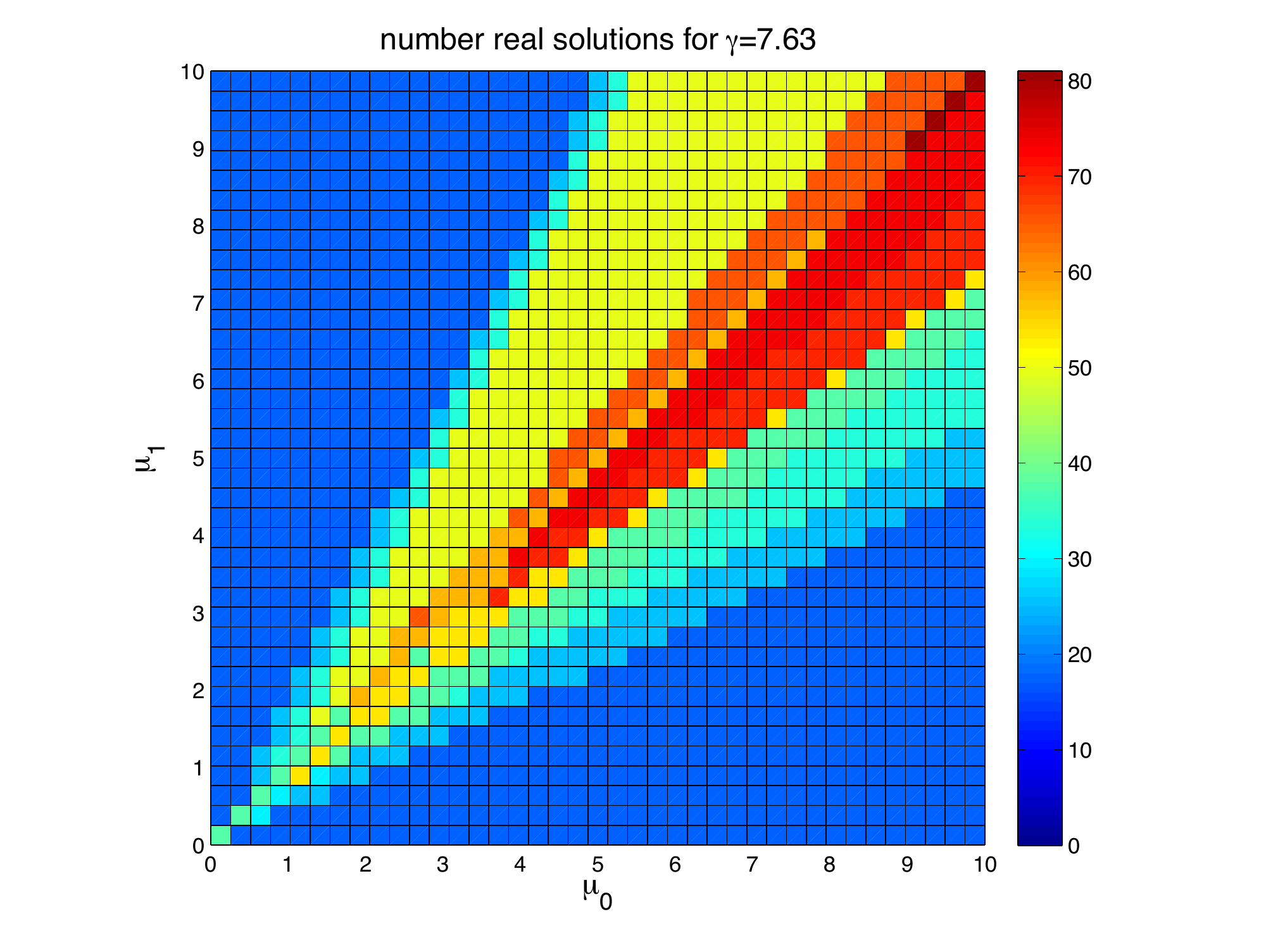}
\caption[Monks Equations (\# solutions) plots]{Regions with various numbers of real solutions to the Monks equations in parameter space ($\mu_0$, $\mu_1$, $\gamma$).  The displayed region is for $(\mu_0\,,\mu_1) \in [0\,,10]\times[0\,,10]$.  At $\gamma = 0 $ there is one solution everywhere, namely the 0 solution; this plot omitted.  For low $\gamma$, there are few solutions near the parameter-origin, as in (a).  As $\gamma$ increases, regions of higher numbers of solutions appear, and these regions move toward the origin. }
\label{fig:monksnumsolns}
\end{center}
\end{figure}

Aside from timing results, Paramotopy is useful in describing the dependence on the number of real equilibria in the system in relation to the values of the parameters, as seen in Figure~\ref{fig:monksnumsolns}.  For generic parameter values, the Monks Equations have 81 distinct isolated complex solutions, with no positive-dimensional components in the variety corresponding to equilibria solution, as verified using a positive dimensional run in Bertini.  Rarely enough, the single variable group total degree start system constructed by Bertini happens to have 81 paths to track as well.  In this figure, we see one slice along $\gamma = 7.63$ in the three dimensional parameter space, with $(\mu_0\,,\mu_1) \in [0\,,10]\times[0\,,10]$. 


%
%

\subsection{Control}
\label{ex:duffing}

The following system is derived from a receding horizon optimization problem, concerned with driving a nonlinear Duffing oscillator to rest \cite{duffing}. The problem depends on parameters, and can be formulated as a polynomial system, for which we desire the roots. 
Paramotopy was designed specifically to solve such a system at multiple parameter points of interest.  We will use one particular version of the system, corresponding to looking two steps ahead.  
The system is given below in Input~\ref{duffing18}.

\begin{figure}[t]
\begin{center}
\subfloat{\includegraphics[width =0.56\linewidth]{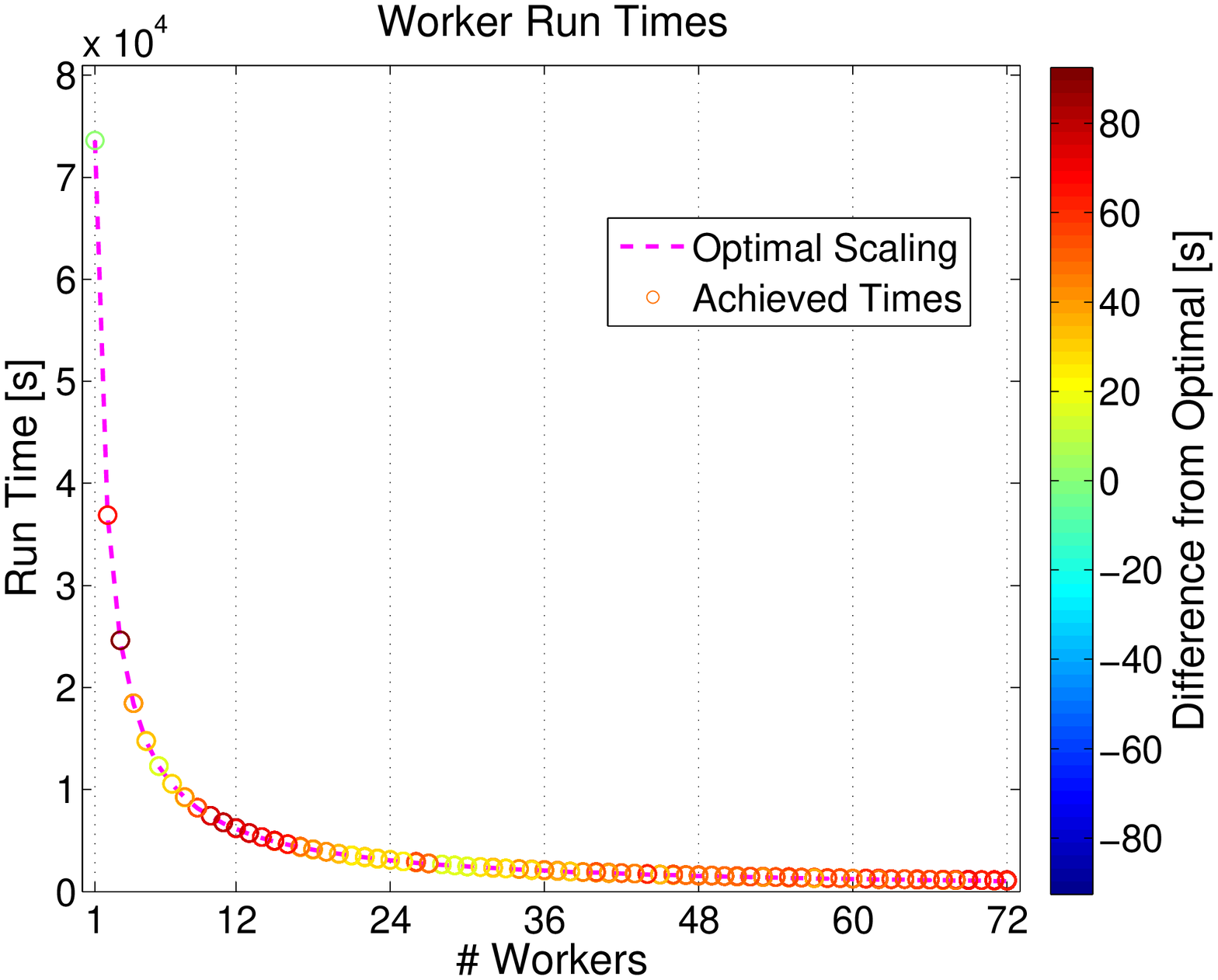}}
\subfloat{\includegraphics[width =0.46\linewidth]{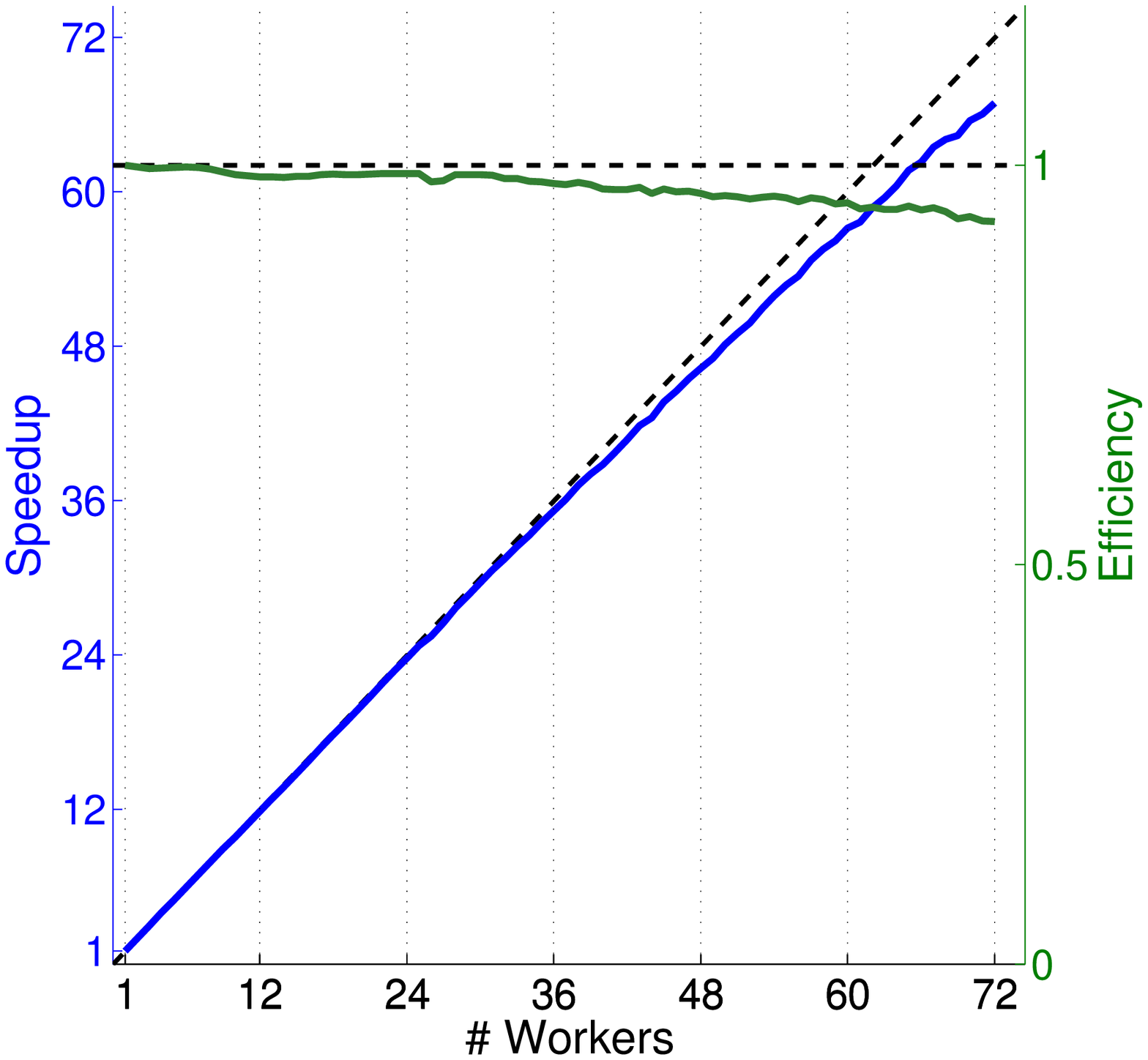}}
\caption[Paramotopy timing -- Duffing $18\times18$ Equations ]{Duffing $18\times18$ timing results.  Linear speedup achieved throughout the range of processors used.}
\label{fig:duffing18timing}
\end{center}
\end{figure}

Timing runs using a mesh of $10^6$ points are demonstrated in Figure~\ref{fig:duffing18timing}.  Again, linear speedup was achieved with efficiency of nearly one throughout the range of processors.  The head spends most of its time waiting for the workers to ask for more work, and the trends in the timing indicate that Paramotopy in this case would have scaled well  beyond 72 workers.

\begin{figure}[t]
\centering
\subfloat[]{\includegraphics[width = 0.49\linewidth]{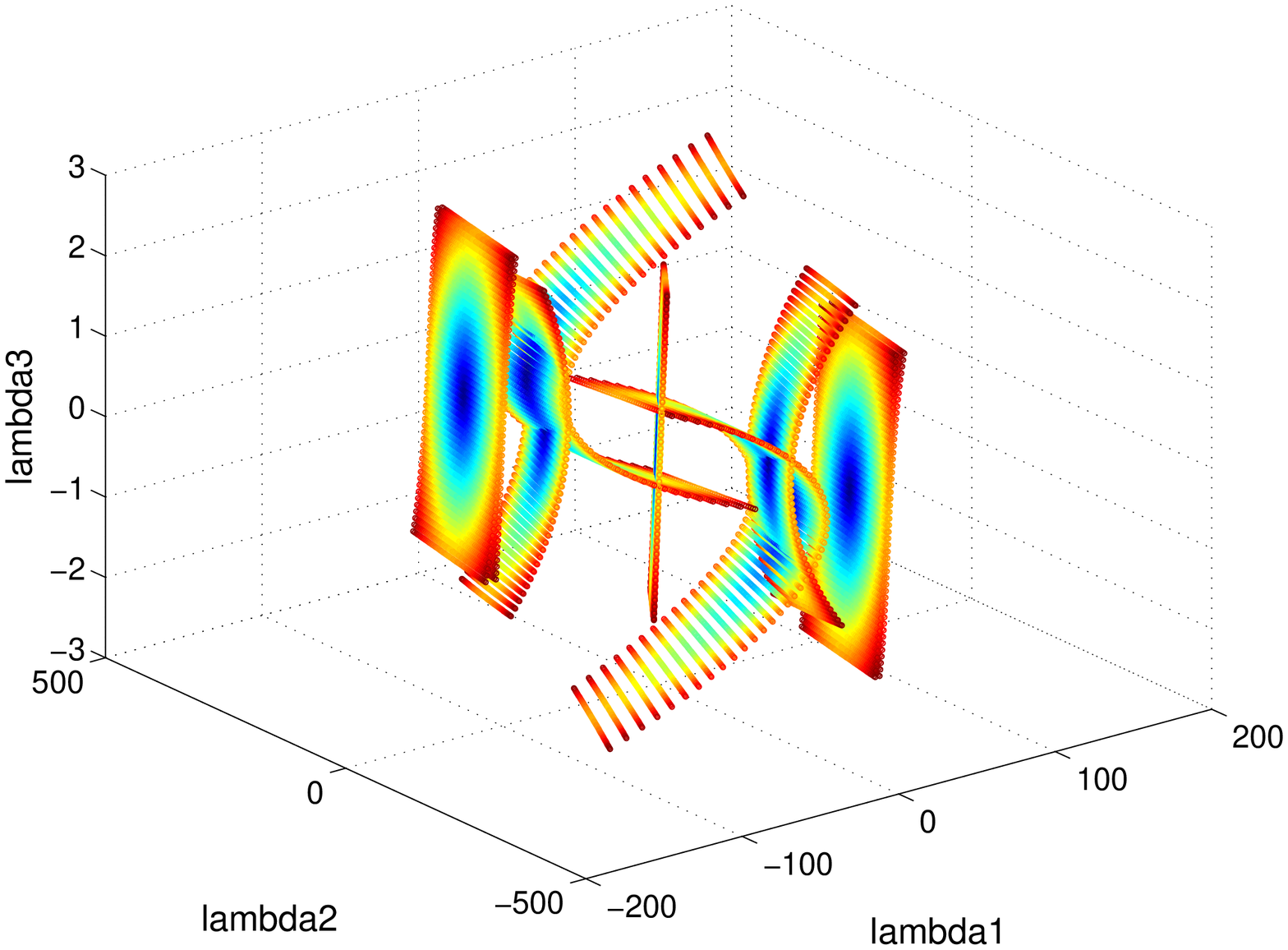}}
\subfloat[]{\includegraphics[width = 0.49\linewidth]{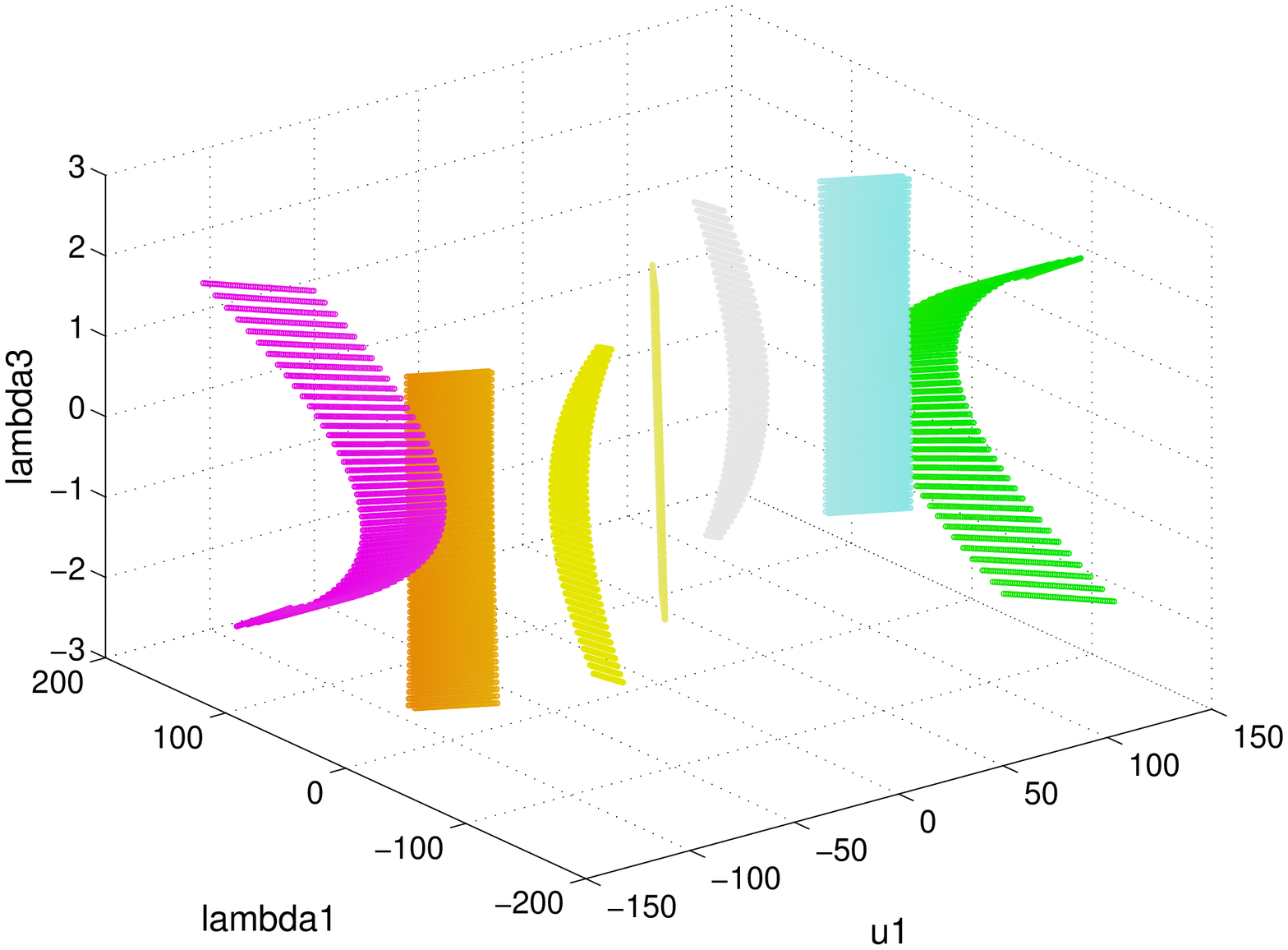}}
\caption[Duffing 18$\times$18 Solution Plots]{Duffing 18$\times$18 solution plots. (a) Projection of solutions onto coordinates $(\texttt{lambda1},\,\texttt{lambda2},\,\texttt{lambda3})$, with color determined by corresponding parameter distance to the origin. (b) Projection onto $(\texttt{lambda1},\,\texttt{u1},\,\texttt{lambda3})$, with RBG value determined by a triple composed of normalized values of $(\texttt{mu3},\,\texttt{mu4}\,,\texttt{x23})$}
\label{fig:duffing18_solutions}
\end{figure}

\subsection{Path failure}
\label{ex:cube}

To demonstrate the capabilities of Paramotopy to deal with path failures, we present here the results of a run using the `cube' system:
\begin{align}
x^6 + y^6 + z^6 -1 = 0 \label{eqn:cube}
\end{align}
Treating $x$, $y$ as parameters, and letting $z$ be the sole variable, we get a one-variable system.  See Input~\ref{cube}.

\begin{figure}[t]
\begin{center}
\includegraphics[width = 0.6\linewidth]{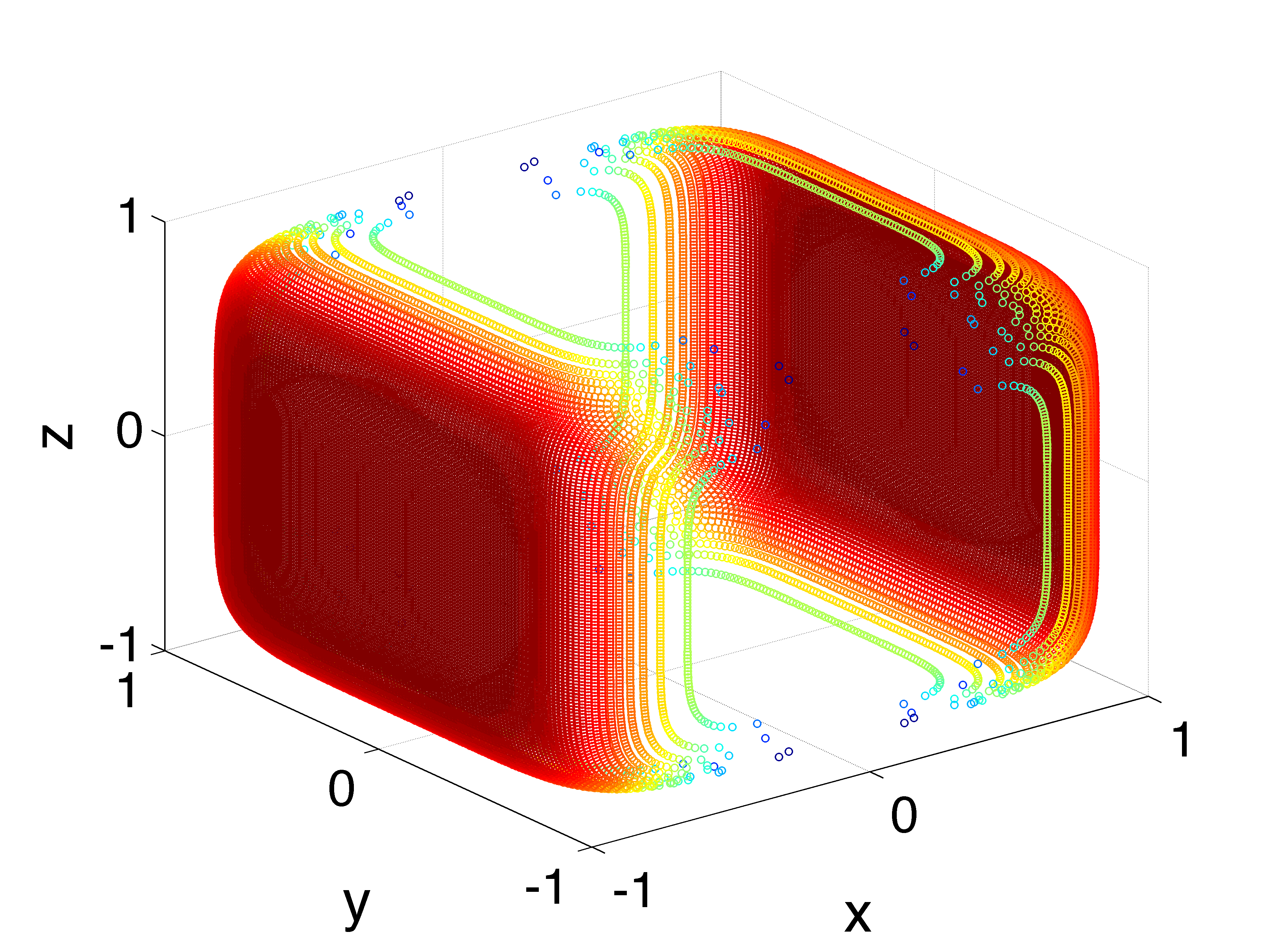}
\caption[Cube]{Plot of samples of $x^6 + y^6 + z^6 - 1=0$.  This system encounters path failures using Bertini default settings, and all these failures can be overcome post-solve using Paramotopy's path failure analysis.}
\label{fig:cube}
\end{center}
\end{figure}

Discretizing $[-1.5,\,1.5]^2$ in a 200x200 grid, for a total of 40,000 points, and solving using the default Bertini configuration, exactly four parameter points result in path failures: 
\begin{align*} (x,y) \in \{ &(0.9425,\,-0.9726),\,      (0.9275,\,-0.8826), \\  &(-1.1536,\,0.7615),\,(-1.1236,\,0.8225) \}.
\end{align*}  Paramotopy re-solved these points, by first moving the random complex start parameter values to another point $p_0$, and then tracking to these four points, using a larger {\tt maxnorm} (the maximum infinity-norm value of a solution at any point during Bertini solve), and {\tt securitylevel 1}.  Of course, whether to move to a new $p_0$ is determined by the user, as are the tolerance values and other Bertini settings.    A plot of the solutions found, concatenated with the parameter values used to give a triple of coordinates, is presented in Figure~\ref{fig:cube}.

One of the benefits of having automatic path failure detection and solution methods is speed.  Solution of the system for all parameter values at the tighter tolerances would have taken much longer; it can be more time efficient to run at lower tolerances for a first pass, and re-solve with more secure settings only when and where needed. 

\section{Conclusions}

In this article, we described the new open-source software package Paramotopy, which can be used to solve 
parameterized polynomial systems very efficiently for large numbers of parameter values.  This extends the reach of 
numerical algebraic geometry in a new direction, particularly a direction that might be useful for mathematicians, 
scientists, and engineers who would like to rapidly test a hypothesis or would like to find regions of a parameter space 
over which the polynomial system has the same number of solutions.  While Bertini and PHCpack have some 
parameter homotopy capabilities, Paramotopy has been optimized for the scenario of using many-processor 
computers to solve at many parameter values of interest.

Paramotopy is under ongoing development, and we expect several extensions in the coming months and years.
For example, Bertini has been improved since the initial development of Paramotopy began in 2010.  In particular, 
Bertini now has configuration settings to turn off parsing (Paramotopy currently uses a special version of Bertini 
built before this new configuration was added), and the parameter homotopy functionality of Bertini has now been 
expanded to include the use of projective space.  

Also, we intend to expand Paramotopy to automatically search for and map out boundaries between cells of the 
parameter space within which the parameterized polynomial system has the same number of real solutions.  More 
precisely, the discriminant locus of a parameterized polynomial system breaks the parameter space into 
cells.  For all points within the same cell, the parameterized polynomial system has the same number of real solutions.  
By taking a grid within some region of the parameter space, successive refinements of this grid near parameter 
values for which solutions are ill-conditioned will help us to ``zoom in'' on portions of the discriminant locus.  This 
``zooming in'' procedure is currently under development and will provide a fast, numerical way of finding these 
boundaries.  Currently, such boundaries can only be computed algebraically (this computation breaks down for 
systems of even moderate size) or through subdivision methods (which have difficulty in zooming in on 
positive-dimensional solution sets such as the discriminant locus).

In the numerical analysis of dynamical systems, significant research has been devoted to the computation of bifurcation points.  MATCONT is a MATLAB package for the bifurcation analysis of ODEs \cite{matcont}.  Another software package for bifurcation analysis is DDE-BIFTOOL \cite{dde}.  Both these packages are similar to some of the capabilities of Paramotopy but both packages are remarkably different from Paramotopy in scope and intended use.  Both of these packages include continuation techniques over $\mathbb{R}$.  Paramotopy is significantly different as it performs continuation methods over 
$\mathbb C$ and works for general polynomial systems.

\section*{Acknowledgments} The authors would like to thank the useful comments from several anonymous referees and Andrew Sommese, 
which have greatly contributed to a better paper.  The first author would also like to recognize the hospitality of 
Institut Mittag-Leffler and the Mathematical Biosciences Institute.


\appendix

\section{Example input files}

This appendix contains the input files for the four examples contained above in Section~\ref{s:ex}.

\DFile{Kinematic system}{robo}{threelinkspatial.txt}
\DFile{Monks system.}{monks}{monks.txt}
\DFile{18$\times$18 Duffing system}{duffing18}{duffing18}
\DFile{Cube system}{cube}{cube.txt}

\bibliographystyle{plain}
\bibliography{paramotopy}

\end{document}